\documentclass[12pts]{amsart}
\usepackage{latexsym,amsthm,amssymb,amsmath,amsfonts,euscript,amstext}
\usepackage{graphicx,epsfig,epstopdf}
\newtheorem{theorem}{Theorem}[section]

\theoremstyle{definition}

\newtheorem{prop}[theorem]{Proposition}

\theoremstyle{remark}
\newtheorem{remark}[theorem]{Remark}
\usepackage{color}

\usepackage[utf8]{inputenc}

\numberwithin{equation}{section}

\date{}

\newcommand{\R}{\mathbb{R}}

\begin{document}
\title{Hartman-Grobman Theorem for Stochastic Dynamical Systems}
\author{Paul Bekima}

\address{ Paul Bekima\newline
Department of Mathematics\\
Morgan State University\\
Baltimore, MD 21251, USA}
\email{pabek1@morgan.edu}

\maketitle

\section*{Introduction}

\vspace{8 pt}

In deterministic systems, the Hartman-Grobman Theorem establishes the "topological equivalence" of the local phase portrait between a system and its linearization around hyperbolic fixed points; simplifying as a consequence the study of the stability of those points. In what follows, we will extend this important theorem on systems perturbed with white noises. Our Stochastic Differential model is of the form \\\

$ dX_{t} = f(X_{t})dt + \epsilon \sigma (X_{t})dB_{t}$ \\

with initial value $ X_{ \{t=0\}} = X_{0} $, for $ 0 \leqslant t \leqslant T < \infty $. \\

Where $ f(x) $ and $ \sigma(x) $ are locally Lipschitz, $ \epsilon $ is a "small" positive parameter.  $ B(t), t\geqslant 0 $,  denote a brownian motions vector of size n, $ \sigma(x) $ is an n by n matrix, $ f(x) $ and $ X(t) $ are also n-dimensional vectors.\\\

However, the conditions for a "useful" linear approximation to a non-linear system are not specific to hyperbolic points, indeed Nils Berglund and Barbara Gentz for example conspicuously used it in their book [3]; particularly during their study of white noise perturbed slow-fast dynamical systems. Yet, since our focus is on understanding the behavior of critical points of a system, we need to make sure that the linear approximation of our perturbed system is equivalent in some sense to the perturbed system of the linear approximation of the corresponding deterministic system. \\\

The paper is organized as follow: we first establish the theorem when the  matrix $ \sigma $ is invertible. We continue by examining the case of non-invertible matrices, then follows the case of non square matrices. We then apply the results to the study of Multi-dimensional slow-fast systems done by Berglund and Gentz [3] by weakening regularity conditions. \\\

 Throughout this paper, unless otherwise specified, $ \lbrace \Omega, \mathcal{F},(\mathcal{F}_{t})_{t \geqslant 0}, P \rbrace $ is a complete probability space with a filtration $ (\mathcal{F}_{t})_{t \geqslant 0} $ satisfying the usual conditions; we will usually say $ \Omega $ as a shorter denomination for  $ \lbrace \Omega, \mathcal{F},(\mathcal{F}_{t})_{t \geqslant 0}, P \rbrace $. All the fixed points are assumed to be hyperbolic. \\

 \newpage

\section{ Hartman - Grobman Stochastic extension }

\vspace{8 pt} 

Assuming that all the coefficients involved are sufficiently regular, we state an extension of the Hartman - Grobman theorem as follows: \\\

\begin{prop}
Around hyperbolic fixed points, the system \\

$ dX_{t} = f(X_{t})dt + \epsilon \sigma (X_{t})dB_{t}$ \\

with initial value $ X_{ \{t=0\}} = X_{0} $, where $ 0 \leqslant t \leqslant T < \infty $, and Lipschitz in both t and x is equivalent to a perturbed version of its deterministic system.\\\

\end{prop}

\begin{proof}

\textbf{Case 1:} $ \sigma $ \textbf{is invertible.} \\

Let $ x_{0} $ be a fixed point, we assume that $ x_{0} = 0 \in \R^{n}, and \ f(x_{0}) = 0 $. \\

The equation becomes \\

\begin{equation}
dX_{t} =( \frac{\partial f}{\partial x}(0) X(t)+ (\frac{\partial^{2} f}{\partial x^{2}})(X^{*}(t), X^{*}(t)) )dt + \epsilon \sigma (X_{t})dB_{t}
\end{equation}

with initial value $ X_{ \{t=0\}} = X_{0} $, where $ 0 \leqslant t \leqslant T < \infty $. \\

Where the convention of repeated indices applies. \\

Let $ L_{t} $ such that $ dL_{t} = \frac{1}{\epsilon} \sigma^{-1} (X_{t}) (\frac{\partial^{2} f}{\partial x^{2}})(X^{*}(t), X^{*}(t))d B_{t} $ \\\

then, \\

$ \langle d X_{t}, d L_{t} \rangle = (\frac{\partial^{2} f}{\partial x^{2}})(X^{*}(t), X^{*}(t))dt $ \\\

Let $ D_{t} = \mathcal{E(L)}_{t} = exp(L_{t} - \frac{1}{2} \langle L_{t}, L_{t} \rangle) $ \\\

By the Girsanov theorem, if $  \mathcal{E(L)}_{t} $ is uniformly integrable, with respect to the probability Q defined by \\

$ dQ = D_{t} dP $, $ X_{t} $ can be written as \\

\begin{equation}
dX_{t} = \frac{\partial f}{\partial x}(0) X(t) dt + \epsilon \sigma (X_{t})dW_{t}
\end{equation}

Where, \\

$ W_{t} = B_{t} - \langle B_{t}, L_{t} \rangle $ \\\

Now it's left to prove that $  \mathcal{E(L)}_{t} $ is uniformly integrable. \\

Since all the functions are continuous, hence bounded on $ [0, T] $ the Novikov's condition $ E [exp \frac{1}{2} \langle L, L \rangle_{T}] < \infty $ is satisfied. \\

\newpage

\textbf{Case 2:} $ \sigma $ \textbf{is non-invertible.} \\

Suppose $ dim ( Im (\sigma(x_{0}))) = k $

Let $ k < n $ and $ v_{1}, v_{2},...v_{k} $ a base of $ Im (\sigma(x_{0})) $, suppose there are k processes $ P_{1}, ..., P_{k} $ such that \\

$ ( \frac{\partial^{2} f}{\partial x^{2}})(X^{*}(t), X^{*}(t)) = \sum_{i=1}^{k} P_{i}(t) v_{i} $. \\

For $ 1 \leqslant i \leqslant k $, let $ a_{i} $ the pre-image of $ v_{i} $ with respect to $ \sigma (x_{0}) $ \\

Letting $ dL_{t} = \frac{1}{\epsilon} \sum_{i=1}^{k} P_{i} dB^{i}_{t} a_{i} $, the same conclusion holds as in the previous case. \\\
 
\end{proof}

\begin{remark} 

The Proposition holds for non-square matrices provided that : \\

In case 1, the invertibility condition be replaced by the condition of full rank, the result follows the use of the pre-image as in case 2. \\\

\end{remark}

\newpage

\section{Multi-Dimensional Slow-Fast Systems}

In chapter 5 of [3] Berglund and Gentz studied the following system:

\begin{equation}
    \begin{cases}
      dx_{t} = \frac{1}{\epsilon}f(x_{t}, y_{t}) dt + \frac{\sigma}{\sqrt{\epsilon}})F(x_{t}, y_{t})dB_{t} \\\

   dy_{t} = g(x_{t}, y_{t}) dt + \sigma^{'} G(x_{t}, y_{t})dB_{t}\\
           
    \end{cases}\,.
 \end{equation} \\

Wher for the sake of uniformity of the notation, we replace $ \lbrace W_{t} \rbrace_{t \geqslant 0} $ with  $ \lbrace B_{t} \rbrace_{t \geqslant 0} $ \\

Here $ \lbrace B_{t} \rbrace_{t \geqslant 0} $ is a k-dimensional Brownian motion. $ \sigma^{'} = \rho \sigma $. \\

The two authors studied a situation where the deterministic system admits a uniformly asymptotic stable slow manifold. They made the following assumptions. \\\

\textbf{Assumptions}

\vspace{8 pt}

\begin{itemize}

\item[1] \textbf{Asymptotically stable slow manifold} \\

\begin{itemize}

\item Domain and differentiability: There are integers  $ n,n,k \geqslant 1 $ such that  $ f \in \mathcal{C}^{2} (\mathcal{D}, \R^{n}) $,  $ g \in \mathcal{C}^{2} (\mathcal{D}, \R^{m}) $,  $ F \in \mathcal{C}^{1} (\mathcal{D}, \R^{n \times\times k}) $,  $ G \in \mathcal{C}^{1} (\mathcal{D}, \R^{m \times k}) $, where $ \mathcal{D} $ is an open subset of $ \R^{n} \times \R^{m} $. We further assume that f,g,F,G and all their partial derivatives up to order 2, respectively 1, are uniformly bounded in norm in $ \mathcal{D} $ by a constant M.\\

\item Slow manifold: There is a connected open subset $ \mathcal{D}_{0} \subset \R^{m} $ and a continuous function $ x^{\star}: \mathcal{D}_{0} \rightarrow \R^{n} $ such that \\

$ \mathcal{M} = \lbrace (x,y) \in \mathcal{D}: x = x^{\star}(y), \ y \in \mathcal{D}_{0} \rbrace $ \\

is a slow manifold of the deterministic system, that is , $ (x^{\star}(y), y) \in \mathcal{D} $ and $ f(x^{\star}(y), y) = 0 $ for all $ y \in \mathcal{D}_{0} $ \\

\end{itemize}

\item[2] \textbf{Non-degeneracy of noise term.} \

The operator norms $ \Vert \bar{X} (y, \epsilon) \Vert $ and $ \Vert \bar{X} (y, \epsilon)^{-1} \Vert $  are uniformly bounded for $ y \in \mathcal{D}_{0} $. \\\

\end{itemize}

\begin{remark}
The domain and differentiability assumption is sufficient to replace the non-linear system by a corresponding linear one. \\\

\end{remark}

Based on those assumptions, Berglund and Gentz stated the main result of that section.\textbf{ See Theorem 5.1.6 (Multidimensional stochastic stable case)}, in page 149 of [3]. \\\

We extend that result by weakening the regularity assumptions as follows: \\\

\begin{prop}
 For  $ n,n,k \geqslant 1 $ such that  $ f \in W^{2, \infty} (\mathcal{D}, \R^{n}) $,  $ g \in W^{2, \infty} (\mathcal{D}, \R^{m}) $,  $ F \in W^{1, \infty} (\mathcal{D}, \R^{n \times\times k}) $,  $ G \in W^{1, \infty} (\mathcal{D}, \R^{m \times k}) $, where $ \mathcal{D} $ is an open subset of $ \R^{n} \times \R^{m} $. \\

 The Theorem 5.1.6 of [3] page 149 still holds provided that the: \\

 Fundamental Condition: $ \int dx E_{x}(1_{A} (X_{t})) \leqslant C \vert A \vert $ is met for a.e.  $ x \in \mathcal{D} $, and $ A \in \mathcal{B}(\R^{n} \times \R^{m}) $. \\
 
 Where $ \mathcal{B }(\R^{n} \times \R^{m}) $ is the Borel set of $ \R^{n} \times \R^{m} $. \\\

  We further assume that regularizations of f,g,F,G and all their partial derivatives up to order 2, respectively 1, are uniformly bounded in norm in $ \mathcal{D} $ by a constant M.\\

\end{prop}

\begin{remark}
The meaning of the Fundamental Condition is that, on set of null measure, the respective modificationx of the coefficients involved do not count a.e. in (x,y). \\
In other words, the laws of the trajectories charge only continuous paths.  \\ 

\end{remark}

\begin{proof}

Let   $ \phi \in \mathcal{C}^{\infty}_{0}(\mathcal{D}) $ be a function, vanishing outside of the unit ball, and such that $ \int_{\mathcal{D}} \phi = 1 $, let's define $ \phi_{\varepsilon} = \frac{1}{\varepsilon^{m}} \phi (\frac{\cdot}{\varepsilon}) $ \\\

Let $ f_{\varepsilon},\ g_{\varepsilon},\ F_{\varepsilon}, \ and \ G_{\varepsilon} $ be the respective regularizations components by components of f, g, F, and G. \\\

The theorem holds with $ f_{\varepsilon},\ g_{\varepsilon},\ F_{\varepsilon}, \ and \ G_{\varepsilon} $. \\

The conclusion follows from the Ascoli-Arzela theorem. \\

Indeed, the system (2.1) can be written as $ dX_{t}  = dA_{t} + d M_{t}$ \\

The challenging part is the $ A_{t} $, however from the Fundamental Condition, \\

 $ \int dx E[\vert \frac{d A_{t}}{dt} \vert^{q}] = C \int[\vert b \vert^{q}] < \infty $. \\

 If we write $ M_{t} $ under the form  $ M_{t} = \int_{0}^{t} \sigma (X_{s})dB_{s} $ \\\

 $ \int dx \Vert M_{t} \Vert_{W^{s,p}}^{p}  = E\int dx \int_{0}^{1} dt \int_{0}^{1} ds \frac{\vert M_{t} - M_{s} \vert^{p}}{\vert t - s \vert^{1+sp}}  $ \\

 $ \frac{1}{p} < s < \frac{1}{2} $, the first inequality is the Sobolev imbedding. \\

 $ \vert M_{t} - M_{s} \vert = \int_{s}^{t} \sigma dW_{s} $ \\

 Hence, $ \leqslant \ C \int dx \int_{0}^{1} dt \int_{0}^{1} ds \frac{1}{\vert t - s \vert^{1+sp} } E ( \int_{s}^{t} \sigma^{2})^{p/2} $\\

Since, $ ( \int_{s}^{t} \sigma^{2})^{p/2} = \int_{s}^{t} \vert \sigma \vert^{p} \vert t - s \vert^{(1- \frac{2}{p})p/2} $, then \\

 $ \int dx \Vert M_{t} \Vert_{W^{s,p}}^{p} \leqslant C \int dx \int_{0}^{1} dt \int_{0}^{1} ds \frac{1}{\vert t - s \vert^{1+sp} }  \vert t - s \vert^{(1- \frac{2}{p})p/2} \int E (\int_{s}^{t} \sigma^{p})dx $ \\

  But, $ \int E (\int_{s}^{t} \sigma^{p})dx \leqslant C \int \vert \sigma \vert^{p} dx \vert t - s \vert $ \\
 
 Hence, $ \leqslant C \int_{0}^{1} dt \int_{0}^{1} ds \frac{\vert t-s \vert^{p/2}}{\vert t-s \vert^{1+sp}} \ \leqslant Constante $ \\

 Hence, the Ascoli-Arzela Theorem can be applied.

\end{proof}

\newpage

\section{Discussion}

The conditions in Proposition 2.2 are not optimal.\\

 Other functional spaces can be explored in particular the space of BV (Bounded Variation) functions. \\

It will also be interesting to deduce results for the system 2.1 from the corresponding "transport system"  which is a second order Partial Differential system, more specifically a Fokker-Planck-Kolmogorov system.

\newpage

\section*{References}

\vspace{8 pt}

\begin{enumerate}

\item[1.]  Strogatz Steven H., Nonlinear Dynamics and Chaos, Third Edition, CRC Press, 2024. \\

\item[2.] Le Gall J.- F., Mouvement Brownien, Martingales et Calcul Stochastique. Springer, 2013.\\

\item[3.] Berglund N, Gentz B., Noise-Induced Phenomena in Slow-Fast Dynamical Sysyems, Springer, 2005. \\

\item[4.] Lions P.L., Equations et Systemes Paraboliques, Lecture at College de France, 2012-2013. \\

\item[5.] Stroock D., Varadhan S., Multidimentional Diffusion Processes, Reprint of the 1997 Edition, Springer, 2006. \\

\item[6.] Mao X, Stochastic Differential Equations and Applications, 2nd Edition, Woodhead Publishing,2011.\\

\item[7.]  Gilbarg D., Trudinger N., Elliptic Partial Differential Equations, 3rd Edition. Springer, 1998.\\

\item[8.]  Revuz D. and Yor M., Continuous Martingales and Brownian Motion, 3rd Edition. Springer, 1999.\\

\end{enumerate}

\end{document}